\documentclass{gtpart}     
\usepackage{pinlabel}  
\usepackage{graphicx}  
\usepackage{amscd}

\usepackage{amsmath}
\usepackage{amssymb}
\usepackage{amsthm}
\usepackage{tabularx}
\usepackage{longtable}
\usepackage{booktabs}
\usepackage{tikz-cd}
\usepackage{subcaption} 

\title{From near-symplectic constructions to trisections of $4$--manifolds}

\author{David T Gay}

\volumenumber{}
\issuenumber{}
\publicationyear{}
\papernumber{}
\startpage{}
\endpage{}
\doi{}
\MR{}
\Zbl{}
\received{}
\revised{}
\accepted{}
\published{}
\publishedonline{}
\proposed{}
\seconded{}
\corresponding{}
\editor{}
\version{}

\theoremstyle{definition}

\newcommand{\id}{\mathop{\rm id}\nolimits}

\begin{document}

\begin{abstract}    
\end{abstract}

\maketitle

\section{Chatty Introduction (which you can skip if you want to jump to the math)}

Rob Kirby was my PhD supervisor at UC Berkeley; I finished my PhD there in 1999. During my first postdoc position, at the University of Arizona, my postdoc advisor Doug Pickrell suggested that I invite Rob to give a colloquium talk, I think that was in the fall of 2001. I remember Rob gave a survey of the beginning of Heegaard Floer theory and in particular I clearly remember him saying something about organizing the various spin$^\mathbb{C}$ structures on a $4$--manifold into ``buckets''. However, the key event of that visit for me was a full moon night hike that Doug organized at Sabino Canyon. (Part of my point here in the story is to make sure to give Doug due credit for his role in all of this. Another point is to emphasize the importance of hiking in Rob's mathematical life.) During the hike, we ended up talking about my thesis work on handle-by-handle symplectic constructions and I remember Rob commenting, ``Oh, so that's what your thesis was about.'' The point is not that he didn't read my thesis but that the essential idea of my thesis was somehow obscured in the whole mechanics of writing, reviewing, graduating, and so forth; I do know that Rob and I both talked about most of the details of my thesis but somehow we had never stepped back and thought about what it was {\em really} about. And once Rob understood what it was {\em really} about, the ideas started flowing.

In this article I want to tell the (mostly mathematical) story of my collaboration with Rob Kirby that started during that night hike and ended up with our ``discovery'' of trisections of $4$--manifolds~\cite{GKTrisections}. I use ``discovery'' cautiously because much of what we discovered was already known and in the end a lot of our work can be thought of as just finding the right way to organize the ideas and the questions. Our work started hovering around various topological constructions related to the vision of extending $4$-dimensional symplectic tools to a larger class of $4$--manifolds using near-symplectic structures (inspired by Cliff Taubes). We then followed what seemed like the natural mathematical trail to the study of broken Lefschetz fibrations (inspired by Auroux, Donaldson and Katzarkov, Usher and Perutz), then to thinking more generally about generic maps from $4$--manifolds to surfaces (inspired by Lekili), which we ended up calling ``Morse $2$--functions''. (A lot of this involved rediscovering for ourselves in our own way of understanding lots of things already known to singularity theorists and others.) From there we more or less stumbled on trisections of $4$--manifolds as a particularly nice outcome of the theory of Morse $2$--functions. At the end I will give a brief survey of some of the developments surrounding trisections that have happened since our first paper, emphasizing my work with Rob and Aaron Abrams~\cite{GroupTrisections} on group trisections. 

I hope to emphasize throughout the powerful influence of Rob's unique way of seeing the mathematical world on my own mathematical development and the naturalness of the mathematical flow of our collaboration. When I was a graduate student and asked Rob for some advice about what to think about mathematically, I don't remember the exact mathematical context, but I remember clearly that his advice was simultaneously vague, unhelpful and deeply meaningful: ``Follow your nose.'' 

\section{Just the math about what trisections are}

\begin{figure}
\centering
\includegraphics[width=.4\textwidth]{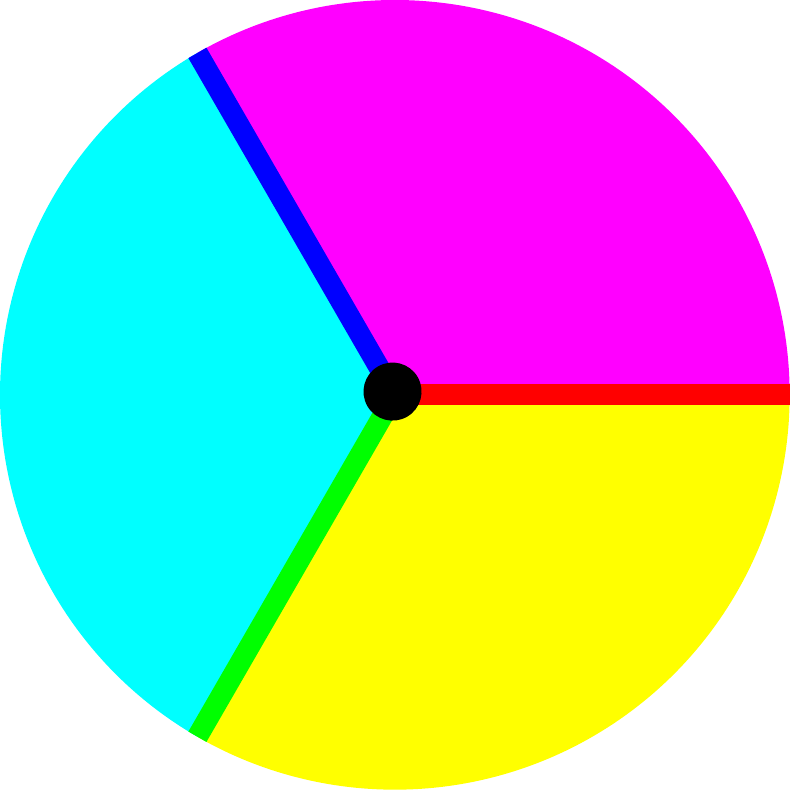}
 \caption{Three slices of a $4$--dimensional pie (pink, yellow and turqoise), their pairwise intersections (red, blue and green), and their triple intersection (black).}
 \label{F:PieSlices}
\end{figure}

For those who have not thought about trisections before, the idea is almost stupidly simple. A trisection of a smooth, closed, oriented, connected $4$--manifold $X$ is a decomposition of $X$ into three topologically simple pieces $X = X_1 \cup X_2 \cup X_3$ which fit together like three ``slices of a pie'' as in Figure~\ref{F:PieSlices}: Each piece $X_i$ is diffeomorphic (after some mild smoothing of corners) to a regular neighborhood of a bouquet of (zero or more) circles in $\mathbb{R}^4$. (A bouquet of zero circles is a point, one circle is a circle, two circles is a figure eight, etc.) Each pairwise intersection $X_i \cap X_j$ is diffeomorphic to a regular neighborhood of a bouquet of (zero or more) circles in $\mathbb{R}^3$  (a.k.a. a solid handlebody). The triple intersection $X_1 \cap X_2 \cap X_3$ is diffeomorphic to a closed, orientable, connected surface. 

\begin{figure}
\centering
\includegraphics{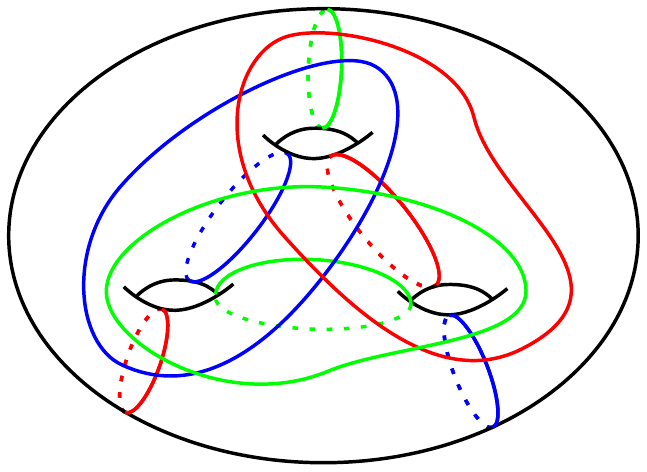}
 \caption{An interesting trisection diagram.}
 \label{F:Genus3Diagram}
\end{figure}

The surprising facts which we discovered between 2011 and 2013 are that {\em every smooth, closed, oriented, connected $4$--manifold} has a trisection ({\em existence}), and in fact it has lots of them; that there is a simple to describe stabilization operation that increases the complexity of a trisection; and that any two trisections of a given $4$--manifold have a common stabilization ({\em uniqueness}). Furthermore, this means that every such $4$--manifold can be described up to diffeomorphism by a {\em trisection diagram} which is a diagram drawn on a closed, oriented surface of genus $g$ involving $g$ red simple closed curves, $g$ blue simple closed curves, and $g$ green simple closed curves; see Figure~\ref{F:Genus3Diagram} for an interesting example. The surface is the triple intersection and the curves of each color are curves that bound disks in the pairwise intersections.

One formally satisfying feature of trisections is that they turn out to be the right generalization of the notion of a Heegaard splitting of a $3$--manifold to the world of $4$--manifolds, with the surprise being that in order to generalize a Heegaard splitting (a decomposition of a $3$--manifold into two simple pieces) one needs to consider decompositions into {\em three} pieces. In fact Heegaard splittings of a $3$--manifold $M$ are best seen as coming from particularly nice smooth maps $f: M \to [-1,1]$ (Morse functions with well organized critical points) and then decomposing $M$ into $f^{-1}([-1,0])$ and $f^{-1}([0,1])$. Similarly, trisections of a $4$--manifold $X$ are best seen as coming from particularly nice smooth maps $f: X \to D^2$ to the disk, called Morse $2$--functions, with particularly well organized singularities, and pulling back the decomposition of the disk into three (fat) pie slices. Also, Heegaard splittings allow $3$--manifolds to be described by {\em Heegaard diagrams}, which are just like trisection diagrams but only involve curves of two colors. Finally, note that a trisection diagram gives three Heegaard diagrams, by considering each pair of two colors.

\begin{figure}
 \centering
 \includegraphics[width=.5\textwidth]{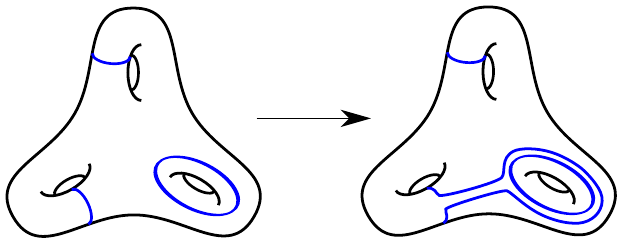}
 \caption{A handle slide on a surface.}
 \label{F:Slide}
\end{figure}

\begin{figure}
\centering
\includegraphics[width=.8\textwidth]{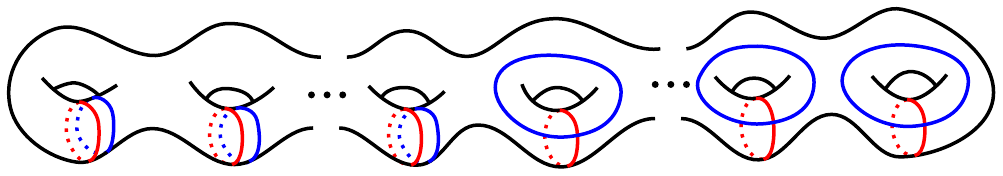}
 \caption{The standard, ``boring'', Heegaard diagram, with some number of parallel red-blue pairs on the left and some number of dual red-blue pairs on the right.}
 \label{F:Standard}
\end{figure}

In principle, a trisection diagram is not much different from a framed link diagram, but it is better organized, and there are algorithms to turn either one into the other. In particular, trisection diagrams have an interesting formal property summarized in the slogan {\em ``pairwise boring, triply interesting''}. There is a standard move on a collection of disjoint simple closed curves on a surface, called a ``handle slide'', which involves one curve ``sliding over'' and picking up a copy of another curve, as shown in Figure~\ref{F:Slide}. Handle slides do not change either the $4$--manifold associated to a trisection diagram or the $3$--manifold associated to a Heegaard diagram. A trisection diagram $(\Sigma, \alpha.\beta,\gamma)$ has the property that each associated Heegaard diagram of two colors $(\Sigma,\alpha,\beta)$, $(\Sigma,\beta,\gamma)$ and $(\Sigma,\gamma,\alpha)$ is equivalent, after some handle slides and a diffeomorphism, to the boring Heegaard diagram in Figure~\ref{F:Standard}, which is a {\em Heegaard diagram} for the connected sum of some number of $S^1 \times S^2$'s (a boring $3$--manifold). However, all $4$--manifolds (and thus a very large class of very interesting objects) can be described by these triples of collections of curves, despite the fact that in pairs they are boring. The interested reader might enjoy verifying this pairwise standardness for the diagram in Figure~\ref{F:Genus3Diagram}.

\begin{figure}
\centering
\includegraphics[width=.2\textwidth]{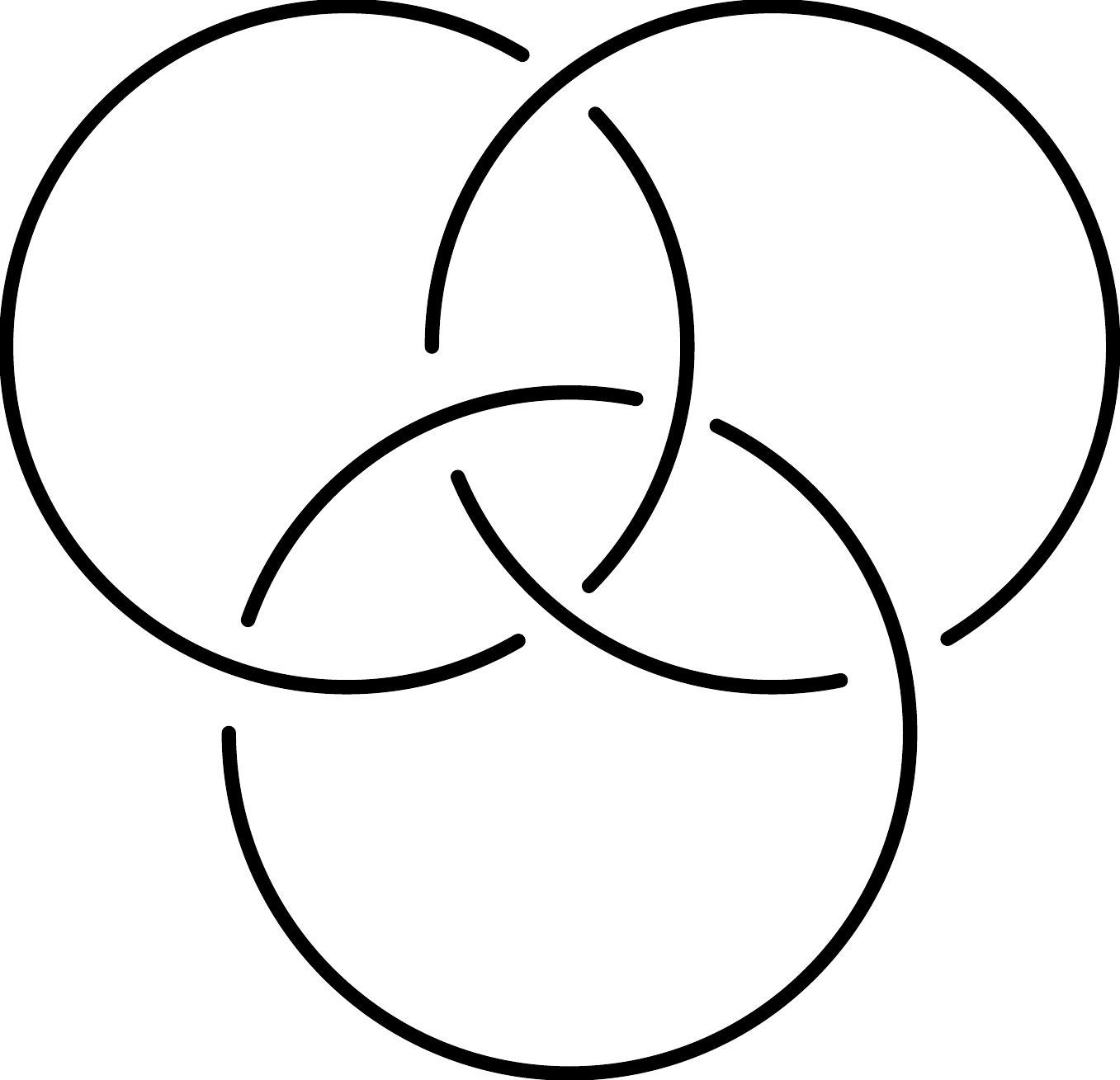}
 \caption{The Borromean rings; each pair of rings is boringly unlinked but the triple is interesting.}
 \label{F:Borromean}
\end{figure}

A classical example of the ``pairwise boring, triply interesing'' slogan, by the way, is the three component link known as the Borromean rings, shown in Figure~\ref{F:Borromean}. This link of three circles has the property than any pair of the three are unlinked but all together they are nontrivially linked. Here is a meta-question for the mathematical community at large: Where else in mathematics does this slogan apply?

\section{The Quixotic quest to make gauge theory mean something topologically}

My entire topological career has been dominated by the incredible twin results of Freedman~\cite{Freedman} that two closed simply connected smooth $4$--manifolds are homeomorphic if and only if they are homotopy equivalent and of Donaldson~\cite{DonaldsonGauge} that there exist invariants that can distinguish pairs (known as ``exotic pairs'') of closed, simply connected smooth $4$--manifolds that are homotopy equivalent but not diffeomorphic. Rob was closely connected with Freedman's work, but Donaldson's work used gauge theory (PDE's coming from physics), and Rob was never a gauge theorist. Ever since I have known Rob an important part of his mathematical quest has been to understand in some ``genuinely topological'' how to see that exotic pairs are not diffeomorphic. The Donaldson invariants were later replaced to some extent with the simpler Seiberg-Witten invariants~\cite{DonaldsonSW} and eventually many gauge theoretic results could be proved using Ozsv\'ath and Szab\'o's Heegaard Floer invariants~\cite{OzsvathSzaboTriples}, but, still, underlying it all is the problem of counting solutions to PDE's. Rob always wanted a way to {\em see} what we were counting explicitly and combinatorially, and to see from a geometric topologist's perspective, in a direct way, why these counts were invariants of smooth structures. 

We would like to be able to stare at some diagrammatic, combinatorial description of a $4$--manifold, ``count'' something in the diagram, perhaps count smooth isotopy classes of some objects in the diagram, and then understand that this is an invariant because one understands the moves needed to relate any two such diagrammatic descriptions of the same $4$--manifold. I would argue that we are still far from such a perspective, but that for formal reasons such a perspective should be out there. The closest we have come that I am aware of are the handful of $4$--dimensional results coming out of Khovanov homology, in particular Rasmussen's computations of slice genus~\cite{Rasmussen} and subsequently Lambert-Cole's trisections and Khovanov homology based proof of the Thom conjecture~\cite{PLCThom}. In what follows, one thread of the story of my collaboration with Rob is that we were constantly trying to catch up with gauge and Floer theory and reinterpret it in ways we non-analytically minded ``mere geometric topologists'' could understand and really see.

Generally, however, almost all invariants proceed by first imposing some additional structure, then counting something, and then showing that the answer does not depend on the choice of that additional structure. In the gauge theoretic context this additional structure is a collection of differential geometric data such as Riemannian metrics, connections on various bundles, and so forth. In the Floer theoretic setting the additional structure is generally an almost complex structure, usually with a symplectic structure lurking in the background. In our hypothetically more ``combinatorial'', more ``topologically understandable'' visions, there would still be extra structure, but this would come in the form of a diagram describing a construction of the manifold or simply a decomposition of the manifold into elementary pieces. For us, inspired heavily by Morse theory, an intermediate structure has typically been some kind of smooth map from the manifold in question to another simpler manifold, e.g. $\mathbb{R}$ or $\mathbb{R}^2$ or $S^2$. Symplectic topology has always lurked in the background as ``borderline understandable'' to simple minded topologists, and thus you will find the following storyline progresses through symplectic topology to understanding various classes of maps from dimension $4$ to dimension $2$ and finally arriving at decompositions of manifolds and associated diagrams, with an addendum in which the ``diagrams'' are eventually replaced with group theory.

There is a great deal also to be said about the extent to which complex algebraic geometry has inspired smooth $4$--dimensional topology, with smooth algebraic surfaces being in some sense the prototypical $4$--manifolds. This is too much for me to competently discuss here but suffice it to say that this thread very much lurks in the background, and I believe that much of Rob's insight into the problems he and I have thought about can be traced back to his early work with Harer and Kas on the smooth topology of complex surfaces~\cite{HarerKasKirby}. (To a complex algebraic geometer, a quick introduction to trisections is to say that they are the natural generalization of the decomposition of the complex projective plane into three coordinate charts.)

Perhaps a better reason for this quest than just the fact that we couldn't understand gauge theory is that all gauge theoretic methods have failed to say anything about homotopy $4$--spheres, and thus shed no light on the holy grail, the smooth $4$--dimensional Poincar\'e conjecture (or its sister, the smooth $4$--dimensional Schoenflies problem). Thus, like any good $4$--dimensional topologist should, we have always held out hope that some day we might stumble upon some invariant, almost certainly not defined using anything like gauge theory, which could distinguish smooth homotopy $4$--spheres and hence disprove the Poincar\'e conjecture. Owning up to such a dream in public is perhaps not good form, but in our collaboration we have had moments of (probably misguided) optimism and enthusiasm in this direction. The truth is that we have been frequently quite naive about the potential to make hard things easy to understand, but that our naivet\'e seems to have served us pretty well.

In what follows I will try to give a fairly thorough account of the path Rob and I took that brought us to thinking about trisections of $4$--manifolds. Of course we never knew where we were going in the long run, but we looked at a sequence of related problems that each raised new questions that, in hindsight, seemed to lead inexorably towards trisections. All of the questions we thought about along the way are, in my opinion, still very important and serve to set trisections in the right context. At each step we were working on $4$--manifolds equipped with certain auxiliary structures, and I will describe those auxiliary structures as they arise.

\section{Near symplectic constructions}

The relevant auxiliary structures to discuss here are symplectic and almost complex structures on $4$--manifolds, building on the foundational work of Gromov~\cite{Gromov}. Symplectic structures are closed, nondegenerate $2$--forms, which means one can always find local coordinates in which they have (in dimension $4$) the form $dx_1 \wedge dy_1 + dx_2 \wedge dy_2$. (A $2$--form at a point is an antisymmetric bilinear form on the tangent space at that point that should be thought of measuring the ``signed area'' of the parallelogram spanned by two vectors, and in this case the area of the $1 \times 1$ squares in the $x_1,y_1$ and $x_2,y_2$ planes are $1$ while the $1 \times 1$ squares in any of the $x_1,x_2$, $x_1,y_2$, $y_1,x_2$ and $y_1,y_2$ planes have ``area'' $0$.) An almost complex structure on a $4$--manifold is a way to ``multiply by $i$'' in the tangent spaces, in other words a smoothly varying linear automorphism $J$ of each tangent space such that $J^2 = -\id$. The structure of a $2$--dimensional complex manifold on a $4$--manifold induces an almost complex structure on the tangent spaces, but not every almost complex structure comes from a complex manifold structure. 

A classical way to probe the structure of a complex algebraic variety is to study enumerative problems for curves in that variety. Remembering that ``curves'' over $\mathbb{C}$ are actually real $2$--dimensional surfaces, one can push these enumerative methods to the softer setting of almost complex manifolds and study ``pseudoholomorphic curves'', which are surfaces in the $4$--manifold whose tangent spaces are fixed by the almost complex structure, in other words, their real $2$--dimensional tangent spaces are actually complex lines with respect to $J$. Gromov's key contribution was to show that, when an almost complex structure is dominated by a symplectic structure (meaning that the $2$--form assigns positive area to a parallelogram spanned by $V$ and $JV$ for any nonzero tangent vector $V$), then one can control families of pseudoholomorphic curves and get compactness results for their moduli spaces. This then allows, in special cases, for one to have well-defined counts of pseudoholomorphic curves and prove that these counts are invariant under various choices. Note that the pseudoholomorphic condition on an embedding of a surface is basically a partial differential equation, so looking for pseudoholomorphic curves is basically counting solutions to a PDE, but somehow it feels a little more concrete than the PDE's involved in gauge theory.

Returning now to the project of finding a more ``concrete'' or ``geometric topologist friendly'' understanding of gauge theoretic invariants of smooth $4$--manifolds, Taubes~\cite{TaubesGWSW} showed that, when a $4$--manifold supports a symplectic structure, its Seiberg-Witten (SW) invariants can be calculated by counting pseudo-holomorphic curves with respect to a generic almost complex structure dominated by that symplectic structure, in other words, a Gromov invariant (Gr). This result is generally summarized as SW=Gr. This has a range of spectacular implications, but when handed a random $4$--manifold, it might not be so clear whether it has a symplectic structure or not. On the other hand, Taubes, Honda~\cite{HondaNearSymplectic} and others had observed that every $4$--manifold with $b_2^+>0$, which means that it contains surfaces which have positive signed intersection number with any wiggled versions of themselves, supports a ``near symplectic'' structure. This is a closed $2$--form which is symplectic away from a $1$--dimensional locus where it is zero, and along this locus it vanishes transversely in the appropriate sense, so that its zero locus is a $1$--manifold in a neighborhood of which the $2$--form has a standard model. Taubes floated the idea that perhaps suitably defined Gromov invariants in this setting would recover Seiberg-Witten invariants for arbitrary $4$--manifolds with $b_2^+$ positive, and initiated a study~\cite{TaubesBeasts} of the behavior of pseudoholomorphic curves in the standard local model near the vanishing locus. This line of reasoning seems to have finally reached fruition with Chris Gerig's work~\cite{Gerig}, which is also a better resource than this paragraph for a proper account of the history and the motivation for the idea.

Still, the fact that every $4$--manifold with positive $b_2^+$ supports a near symplectic structure was not constructive, however, so if the SW=Gr plan panned out in the near symplectic setting it was still not clear how useful that would be if one did not know how to explicitly construct a near symplectic form on a particular $4$--manifold. After our Sabino Canyon hike, Rob explained these ideas to me and then Rob and I set out to solve this problem~\cite{GayKirbyNearSymplectic}, where ``explicitly'' meant to us ``starting from a framed link diagram describing a handle decomposition of the $4$--manifold'' and then proceeding in something that might loosely be called an algorithm. 

It is important to point out here that in the end it is a bit disingenuous to describe our construction as explicit because of one key point; at a critical stage in our construction we had two near-symplectic symplectic structures, one on one half of the $4$--manifold and one on the other, and we needed to glue them together along the separating $3$--manifold. There is a standard way to do this, using a {\em contact structure} on the $3$--manifold as the appropriate gluing boundary data. (If you don't know what a contact structure is, just know that it is the correct boundary data for symplectic structures.) In order to glue, we needed to know that the contact structures coming from the constructions on the two halves were equal (or could be deformed to be equal). We knew this thanks to Eliashberg's result~\cite{EliashbergClassification} that a certain class of contact structures, known as {\em overtwisted contact structures}, could be classified using algebraic topological invariants, so we adjusted our construction appropriately to make sure that the invariants coming from the two sides matched. However it is not at all clear that Eliashberg's machinery to go from knowing that the invariants match to deforming (isotoping) the contact structures to be equal is ``explicit''.

\section{Broken Lefschetz fibrations}

The auxiliary structures of concern here are Lefschetz fibrations; for a good overview see Gompf's exposition in the Notices of the AMS~\cite{GompfWhatIsLF}. In short, a Lefschetz fibration on an oriented $4$--manifold $X$ is a smooth map from $X$ to the $2$--sphere $S^2$ which has finitely many isolated singularities, each of which is locally modelled (respecting orientations) on the simplest singularity that can arise in a holomorphic map from $\mathbb{C}^2$ to $\mathbb{C}$: $(z_1,z_2) \mapsto z_1^2+z_2^2$. Away from these singularities, the map thus looks like a surface bundle over a surface, so locally like $\Sigma \times B^2 \to B^2$ for some surface $\Sigma$. To get a more topological understanding of what it means to support a sympectic structure, Donaldson~\cite{DonaldsonLP} showed that, after blowing up enough times, all symplectic $4$--manifolds support Lefschetz fibrations. (``Blowing up'', despite its name, is a mild operation that changes a manifold in a very controlled way so as to allow two surfaces that intersect to become disjoint, at the expense of introducing a new surface, called the ``exceptional divisor'', that they both intersect; as long as one keeps track of the exceptional divisor one still has access to the original manifold and all of its topology.) With this in mind, Usher~\cite{UsherDSGr} showed that one could recover the full Gromov pseudoholomorphic curve count by a certain count of pseudoholomorphic ``multi-sections'' of symplectic Lefschetz fibrations, arguably making Taubes's SW=Gr result slightly more meaningful to the average geometric topologist, and giving us a  slightly more topological way to think of what the Seiberg-Witten invariants are counting. 

Inspired by Taubes's near-symplectic vision discussed in the preceding section, Auroux, Donaldson and Katzarkov~\cite{ADK} generalized Donaldson's Lefschetz fibration result to show that every near-symplectic $4$--manifold, after blowing up, has the structure of a ``broken Lefschetz fibration'' (BLF), which is like a Lefschetz fibration over $S^2$ but also allows for a $1$--dimensional locus of singularities of a particular model, called {\em indefinite folds}. This $1$--dimensional singular set is exactly the same zero locus where the near symplectic form vanishes. Following our noses, we naturally wondered next how to construct these BLF's ``explicitly'' from, for example, a handle decomposition of a $4$--manifold. In the end~\cite{GayKirbyBALFs} we discovered that if we allowed both orientation preserving and orientation reversing local models for the Lefschetz type singularities (we called the resulting fibrations ``broken achiral Lefschetz fibrations'', or BALFs), then we could construct BALFs on {\em all} closed oriented $4$--manifolds. 

The construction was similar in spirit to our construction of near-symplectic forms: First we found a good way to split the given $4$--manifold into two pieces on each of which we could construct partial fibrations and then we figured out the appropriate boundary conditions to govern these fibrations on the separating $3$--manifold. (Experienced low-dimensional topologists will not be surprised that this boundary condition is the structure of an ``open book decomposition''.) Then we massaged our constructions carefully so as to arrange that this boundary data agreed and we coud glue the fibrations together. As before, it would be disengenuous to describe our construction as explicit because of this last step, which passed from open book decompositions to contact structures by way of the Giroux correspondence~\cite{Giroux} and then appealed again to Eliasherg's classification~\cite{EliashbergClassification} of overtwisted contact structures.

Aside from this caveat about explicitness or the lack thereof, another thread the reader may be picking up is the increasing mixing of categories here: The study of almost complex structures, symplectic structures and Lefschetz fibrations already lies (a little bit uncomfortably) somewhere between complex algebraic geometry and differential topology. Throwing in circles along which new types of degenerations arise, and in the fibration case having both isolated complex type singularities as well as $1$--dimensional decidedly nonholomorphic singularities, and then allowing the isolated singularities to flip orientations, really starts to seem like a bit of a stretch and an awkward mix of perspectives. The next phase in our collaboration started to clean this up.

\section{Morse $2$--functions}

The auxiliary structures of relevance here are, in a broad sense, stable smooth maps from manifolds of some dimension to other manifolds of lower dimension. A map is stable if small perturbations do not change its essential qualitative features; more precisely, the original map and the perturbed are equal after pre- and post-composing with isotopies of the domain and range. Thus the function 
$y=x^2$ is stable in the world of smooth maps from $\mathbb{R}$ to $\mathbb{R}$ while the function $y=x^3$ is not. Stable maps from smooth manifolds to $\mathbb{R}$ (and sometimes to other $1$--manifolds such as the circle or a closed interval) are called {\em Morse functions} and have been tremendously useful tools for probing the topology of smooth manifolds in general, especially when used in conjunction with gradient-like vector fields to give handle decompositions. Stable maps to dimension $2$ are also fairly well understood but have not been used as extensively as Morse functions as a tool to probe the topology of manifolds.

Perutz~\cite{PerutzMatching} studied the problem of defining, in the broken Lefschetz fibration setting, something analogous to the counts of pseudoholomorphic multisections for Lefschetz fibrations studied by Usher, motivated by the possibility that these would then again recover Seiberg-Witten invariants, or at least be invariants even if one did not know that they agreed with pre-existing invariants. To show directly that these counts, called Lagrangian matching invariants, are in fact $4$--manifold invariants, and did not depend upon the choice of BLF, one would need to understand how to move from one BLF on a given $4$--manifold to another. In other words, one needs not only existence results for B(A)LFs, but also {\em uniqueness} results. 

While we were thinking about this (working together at the African Institute of Mathematical Sciences in Muizenberg, South Africa, probably in 2007) we found out about Lekili's insightful observation~\cite{Lekili} that if one wanted to connect two BLF's by a $1$--parameter family of smooth functions, one should first observe that Lefschetz singularities are actually not stable in the world of smooth maps. A small perturbation in the world of smooth maps will make a Lefschetz singularity into a much more complicated circle of singularities. Thus a generic $1$--parameter family of smooth functions connecting two Lefschetz fibrations or BALFs should be expected to pass through functions that are much more general and do not have any Lefschetz singularities at all. The local behavior of stable smooth maps from $4$--manifolds to $2$--manifolds has in fact been understood for a long time and Lekili's paper~\cite{Lekili} has an excellent appendix that runs through the general machinery for understanding stable maps between various dimensions~\cite{Wassermann} applied in the case of dimensions $4$ and $2$. At least this is where I learned this material, Rob probably basically had already absorbed most of this by osmosis over the years and just needed a little reminder from Lekili's article.

Lekili essentially answered the question of how to think about local moves connecting one BALF to another, but the key problem was that in the intermediate stages one might wander quite far away from the world of BALFs and move through maps that are nothing like small perturbations of BALFs. In particular, one might run into what are called {\em definite folds}, discussed in more detail below. Whether one could avoid definite folds is analogous to the question of whether, in a $1$--parameter family of functions connecting two given Morse functions with the same number of minima and maxima, one can avoid introducing extra minima or maxima along the way. In general this can be done, modulo some obvious counterexamples and low-dimensional exceptions. Because of the analogy with Morse theory, Rob and I began calling stable maps to dimension $2$ ``Morse $2$--functions'', thinking of them as vaguely like a $2$--category version of Morse functions, whatever that might mean. (Here we were very much inspired by conversations with Peter Teichner at MSRI in the spring of 2010.) 

The key idea of a Morse $2$--function is that locally a Morse $2$--function $F$ on an $n$--manifold looks like a generic $1$--parameter family $f_t$ of Morse functions on an $(n-1)$--manifold, so that there are local coordinates so that $F(t,p) = (t,f_t(p))$ where $t$ is a single time parameter, $p$ is an $(n-1)$--dimensional ``spatial'' coordinate, and $f_t$ is a generic path of functions connecting two Morse functions. However, globally there is no well-defined time direction, either in domain or range. Generic $1$--parameter families of Morse functions will be honest Morse functions for all but finitely many times, and at finitely many times will experience births or deaths of pairs of critical points, or the coincidence of two critical points having the same critical value at an instant. The tracks of the critical points for the honest Morse functions are called ``folds'', and are $1$--dimensional in domain and range, while the birth/death events are called ``cusps''; the coincidences of having two critical points with the same critical value are folds whose images in the range cross. The generic behavior of homotopies between Morse $2$--functions is exactly locally modelled on the generic behavior of homotopies between homotopies between Morse functions, in other words the moves of ``Cerf theory''.

\begin{figure}
 \centering
 \includegraphics[width=.5\linewidth]{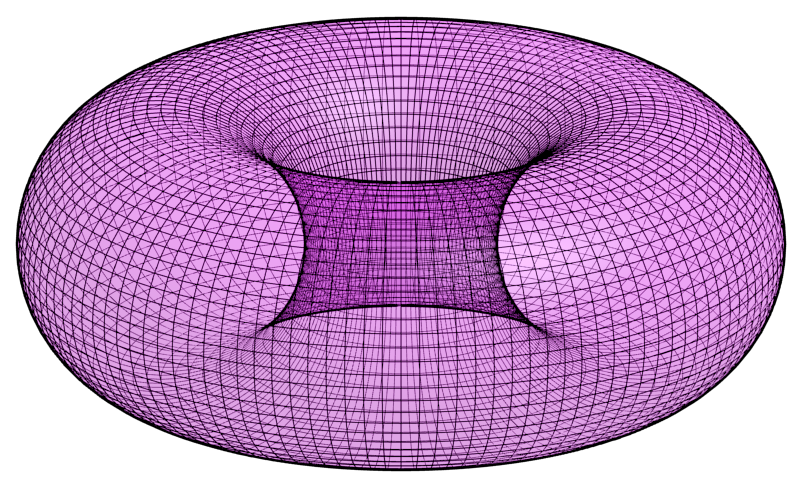}
 \caption{A Morse $2$--function on a torus.}
 \label{F:2DMorse2Fcn}
\end{figure}

If the preceding paragraph did not mean much to the reader, the basic idea is well illustrated with a picture of a Morse $2$--function on a surface. Figure~\ref{F:2DMorse2Fcn} shows a Morse $2$--function on a torus. In fact, any time one draws a picture of a surface on a piece of paper, one is of course presenting a map from that surface to $\mathbb{R}^2$ and, assuming genericity, this will be a Morse $2$--function. The folds are literally the places where the surface folds over, and the cusps are where folds ``switch directions'' in some sense; these are all clearly visible in this illustration. The preimage of a nonsingular point is an even number of points in the torus, and the number of these points jumps by $2$ when crossing a fold. The only difference in higher dimensions in that different dimensions and codimensions can fold in opposite directions along a fold, and the preimages of nonsingular points are higher dimensionsal submanifolds, not collections of points. In dimension $4$, the preimages of points are surfaces, and the topology of these surfaces changes as one crosses a fold.

Note that the index of a fold ranges from $0$ to $n-1$, but is only well-defined up to switching $k$ with $n-1-k$. A fold of index $0$ (equivalently $n-1$) is called a ``definite fold'' and all other folds are called indefinite. Thus on $4$--manifolds there are only really two types of folds: definite index $0$ (equivalently $3$) folds and indefinite index $1$ (equivalently $2$) folds. (On a $2$--manifold, all folds are definite.) One reason definite folds in dimension $4$ might be undesirable is that they in general lead to disconnected fibers, since crossing a definite fold in the index $0$ direction creates a new $S^2$--component of the fiber. This (and $S^2$--fibers in general perhaps) is bad from a symplectic geometry and invariant constructing perspective, as explained to us by Katrin Wehrheim, another major motivator for us in our early work in this subject.

The upshot is that we showed~\cite{GayKirbyMorse2Fcns} how to construct indefinite, fiber-connected Morse $2$--functions on a given $n$--manifold in a given homotopy class of maps, when $n>3$ and a natural $\pi_1$ condition is satisfied, and we also showed that any two such Morse $2$--functions can be connected by a generic homotopy maintaining the indefinite and fiber-connected properties, again when $n>3$. Of course, all along we were primarily interested in the case $n=4$ and, having settled this question, we now started to think seriously about what to do with Morse $2$--functions in general on smooth $4$--manifolds, how to put them into particularly nice forms, and how to use them to appropriately probe the topology of the manifolds involved.

\section{Trisections}

I started my position at the University of Georgia (UGA) in August 2011. One thing Rob was always fantastic about was travelling to visit me wherever I ended up, understanding that my family life with three kids made it harder for me to up and travel at the drop of a hat. And Rob had a long time connection to UGA going back to the early days of the Georgia Topology Conference, so he started visiting me at UGA pretty regularly (and making sure to combine each such visit with some social time with old retired friends Jim Cantrell and John Hollingsworth). Over the course of a few visits in late 2011 and early 2012 we were busy drawing pictures of $\mathbb{R}^2$--valued Morse $2$--functions on homotopy $4$--spheres arising from balanced presentations of the trivial group not known to be Andrews-Curtis trivializable. Our grand hope of course was to find an invariant that could show that these homotopy $4$--spheres were not diffeomorphic to $S^4$, and we had an idea that would involve counting certain surfaces with boundary on the folds. (We called these surfaces ``flow charts'' because they would be swept out by flow lines of gradient-like vector fields in local charts; this is a loose end idea that never panned out for us, just in case the story was seeming too tidy with suspiciously few false starts and dead ends.)

In the process we developed a very systematic and standardized way to organize the image of a $\mathbb{R}^2$--valued Morse $2$--function and, suddenly, one day trisections jumped out at us! In hindsight so much of it is obvious but at the time it seemed somehow  unusual and dramatic. The key point is that we were busy building Morse $2$--functions ``from left to right'' starting from handle decompositions of $4$--manifolds, with the $0$-- and the $1$--handles at the left contributing a very standard part of the Morse $2$--function, the $3$-- and $4$--handles at the right contributing a mirrored standard piece, and the $2$--handles in the middle doing something interesting. Then suddenly trisections jumped out at us when we realized one day in my office at UGA that we could split things a little differently by separating the cusps into three groups.

\begin{figure}
 \begin{center}
  \begin{subfigure}{.4\textwidth}
   \centering
   \includegraphics[width=.8\linewidth]{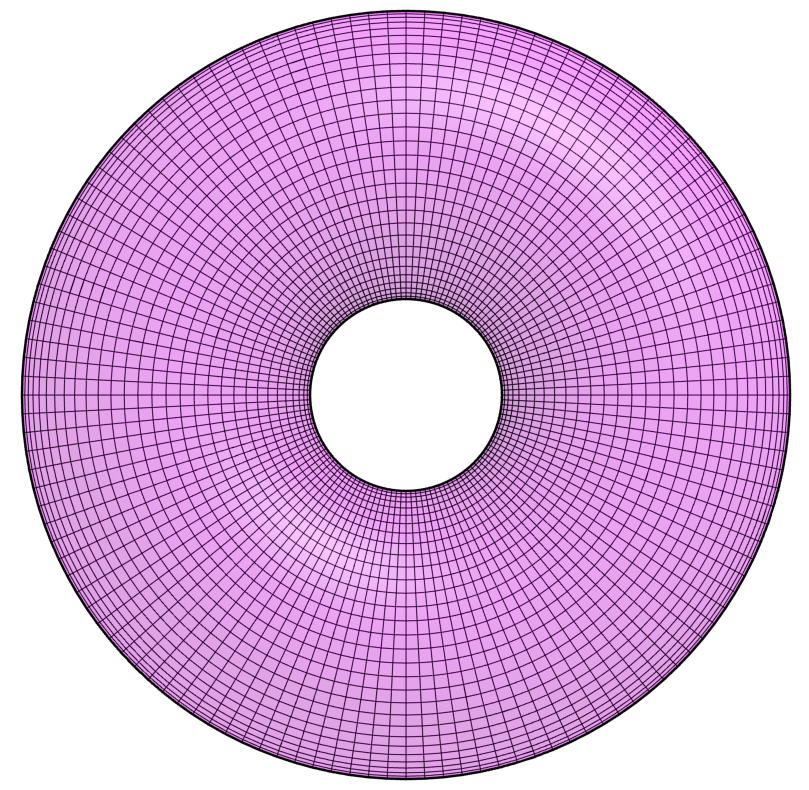}
  \end{subfigure}
  \begin{subfigure}{.4\textwidth}
   \centering
   \includegraphics[width=.8\linewidth]{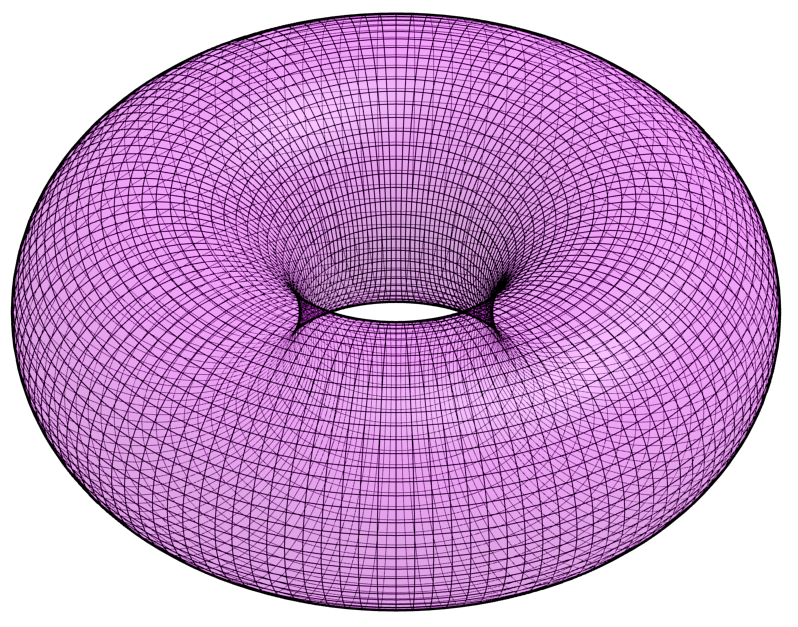}
  \end{subfigure}\\
  \begin{subfigure}{.4\textwidth}
   \centering
   \includegraphics[width=.8\linewidth]{T2V3.png}
  \end{subfigure}
  \begin{subfigure}{.5\textwidth}
   \centering
   \includegraphics[width=.8\linewidth]{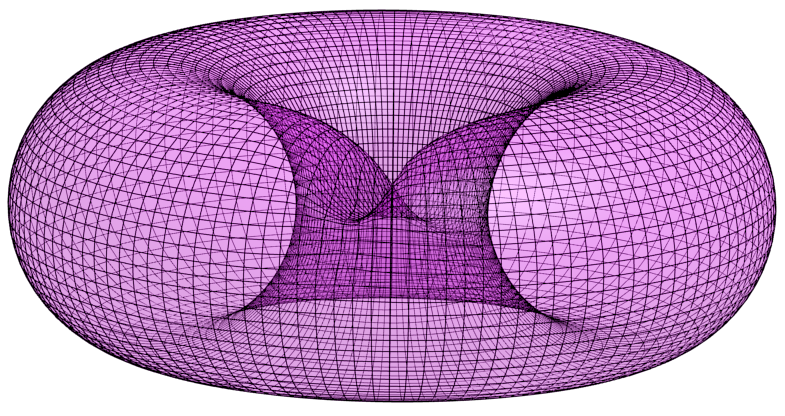}
  \end{subfigure}

 \end{center}

 \caption{A sequence of Morse $2$--functions on the torus.}
 \label{F:SettingUpT2}
\end{figure}

\begin{figure}
 \centering
   \includegraphics[width=.6\linewidth]{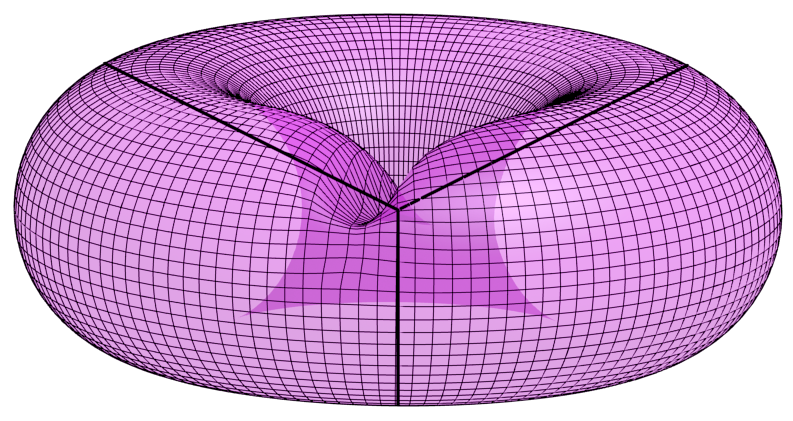}
 \caption{The last Morse $2$--function from Figure~\ref{F:SettingUpT2} with three cut lines indicated in black.}
 \label{F:TrisectingT2CutLines}
\end{figure}

\begin{figure}
   \centering
   \includegraphics[width=.8\linewidth]{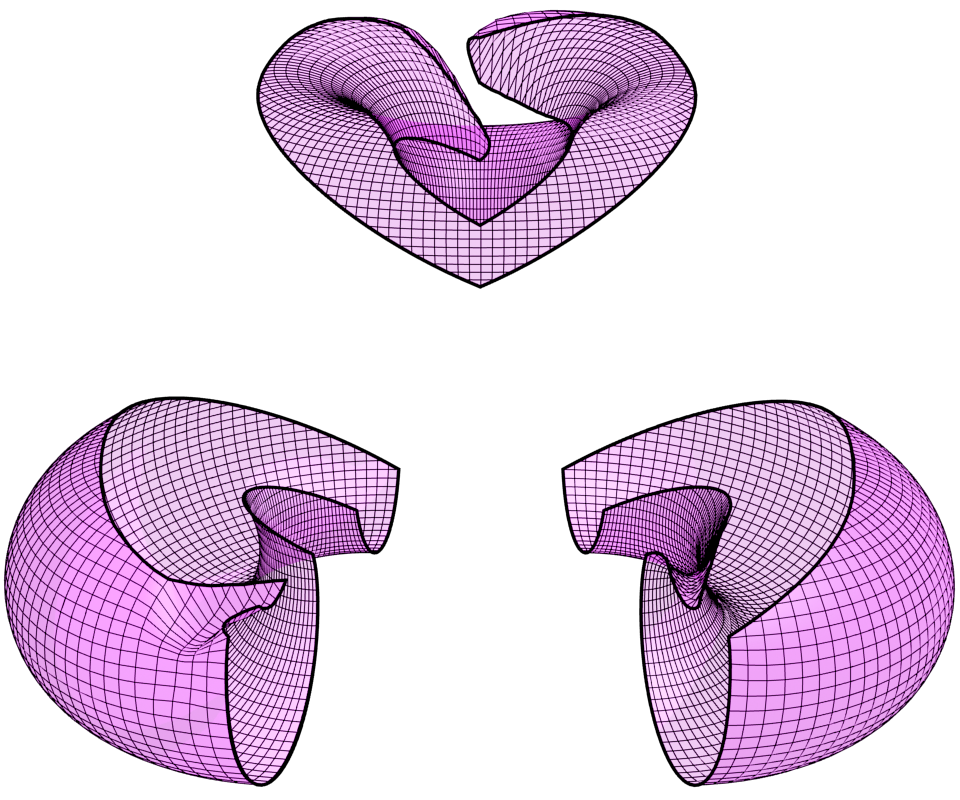}

 \caption{A trisection of the torus coming from cutting open along the preimages of the three black lines in Figure~\ref{F:TrisectingT2CutLines}.}
 \label{F:TrisectingT2}
\end{figure}

The basic idea is well illustrated by a half-dimensional cartoon. (I believe that Rob taught all of his students and collaborators the value of the half-dimensional analogy when working in dimension $4$.) In Figure~\ref{F:SettingUpT2}, we show a sequence of Morse $2$--functions on the torus, each obtained by performing some procedure to the preceding one, and ending with one which might seem unnecessarily complicated. However, if one then trisects the torus using the last Morse $2$--function, as illustrated in Figure~\ref{F:TrisectingT2CutLines} and Figure~\ref{F:TrisectingT2}, we see that in fact the torus is decomposed into three disks (careful inspection shows that these disks are hexagons), and this is the basic ``trisection'' of $T^2$. We essentially mimic this procedure with more care in dimension $4$ to produce trisections of arbitrary $4$--manifolds.

We released two versions of our paper on the arXiv with quite a gap between them and a significant difference. One of Rob's life lessons for me was that if you prove an existence theorem then you should prove a uniqueness theorem. Our first version (May 2012) only had the existence result and then it took us more than a year to nail down the uniqueness result (September 2013). Somewhere in between those two dates was a hardworking session at the Max Planck Institute for Mathematics in Bonn, in which Andras Stipsicz had several helpful comments about how to proceed with uniqueness. Once we had the uniqueness result it was stupidly obvious, and the (for us) really startling fact was just how closely it mirrored the uniqueness result for Heegaard splittings of $3$--manifolds.
%

Having ``discovered'' trisections, it did also become clear to us that in many ways we were obviously rediscovering things others already knew quite well. Of particular note is the fact that Birman and Craggs~\cite{BirmanCraggs} essentially already thought about trisections of $4$--manifolds from a mapping class groups of surfaces perspective, using the term ``triadic $4$--manifolds''. (At least, this term appears in a section heading in that paper, although it never appears in the actual textual body of the paper. Perhaps they had an intuition that there was something more going on with these structures on $4$--manifolds or even, as is so often the case, knew much more than they wrote in their paper.) More recently, Oszv\'ath and Szab\'o were very obviously working with these structures in the guise of Heegaard triples~\cite{OzsvathSzaboTriples}. Our main existence result is more or less implicit in Oszv\'ath and Szab\'o's work, and one could interpret the thrust of our work as being to understand the raw {\em topological} implications of the fact that every $4$--manifold can be described by a Heegaard triple, rather than simply viewing Heegaard triples as a tool to understanding various maps in Heegaard-Floer theory. However the uniqueness result, and the study of trisections as actual decompositions of $4$--manifolds, unifies the subject and sets the stage for many interesting directions of study.

\section{Group trisections}

Many other mathematicians have made striking contributions to the theory of trisections since our first paper. Each of this contributions is worth a lengthy discussion, but since we are focusing on Rob Kirby's work here, I'll briefly describe the most important extension to the theory that Rob and I were involved in, namely the development of ``group trisections''~\cite{GroupTrisections} in collaboration with Aaron Abrams.

One of the first things we thought about with trisections was understanding how much information about a smooth $4$--manifold one retains when applying algebraic topological functors to a trisection, (to all the pieces and their intersections and the various inclusion maps). Since Freedman taught us that in some special cases the algebraic topology of $4$--manifolds seems to basically record homeomorphism type, we expected somehow that by passing from smooth trisections to the algebraic topological shadows of trisections we would pass from smooth information to topological information and then be able to understand something about how trisections of homeomorphic but not diffeomorphic manifolds are related. In fact it turned out, after Aaron Abrams reminded Rob and me of some basic facts from $3$--manifold topology, that the most basic algebraic topology functor of all, $\pi_1$, loses {\em no information at all}! The inclusion maps between the triple intersections, double intersections, single pieces and total manifold making up a trisection fit into a cube of spaces and, when the $\pi_1$ functor is applied, one gets a cube of groups involving one surface group, six free groups, and the fundamental group of the total $4$--manifold and from this cube of groups one can completely recover the smooth $4$--manifold and its trisection. 

We thus defined a ``group trisection'' of a given group $G$ as a cube of groups and homomorphisms between them such that, at one vertex we have the fundamental group of a closed, oriented surface, at the opposite vertex we have the group $G$, and at the 6 intermediate vertices we have free groups. The maps all flow from the surface group in the direction of $G$, all the maps are surjective, and all the six faces of the cube are pushouts. The main result is that a trisected group gives a trisected $4$--manifold and vice versa. This is striking; for example, a trisected group carries with it an integral quadratic form, namely the intersection form of the associated $4$--manifold.

The crucial piece of $3$--manifold topology that Aaron reminded of us is the following ``standard fact'' (for a proof see~\cite{LeiningerReid}): Any surjective group homomorphism from the fundamental group of a genus $g$ surface to a free group of rank $g$ arises as the map induced by the inclusion of the surface as the boundary of a genus $g$ handlebody, and the handlebody filling of that surface is uniquely determined by the group surjection. This, together with Laudenbach and Po\'enaru's result~\cite{LaudPoen} that every diffeomorphism of a connected sum of $S^1 \times S^2$'s extends across a boundary connected sum $S^1 \times B^3$'s is really all you need to prove that a trisected $4$--manifold can be completely recovered from the associated group trisection.

In fact, everything we did is foreshadowed in, and inspired by, John Stallings's famous paper~\cite{StallingsHowNot} ``How not to prove the Poincar\'e conjecture'', which begins with the memorable line ``I have committed the sin of falsely proving Poincar\'e's Conjecture. But that was in another country and until now, no one has known about it.'' One consequence of the dictionary between trisections of groups and trisections of $4$--manifolds is that the smooth $4$--dimensional Poincar\'e conjecture can be reformulated entirely group theoretically, just as Stallings reformulated the $3$--dimensional Poincar\'e conjecture purely group theoretically. I argued unsuccessfully with my coauthors in favor of titling our group trisections paper ``How not to prove the smooth $4$--dimensional Poincar\'e conjecture''. I suspect that Rob, in particular, rejected my proposal because, deep down, he is a fundamentally optimistic mathematician and, whether the conjecture is true or false, he's always willing to consider the possibility that some new idea really might be the right way forward to settling the question once and for all.

\section*{About the author}

David Gay was born in 1968 in Suacoco County in upcountry Liberia, grew up in Liberia, England, Lesotho and Massachusetts, ended up with a PhD in Mathematics from UC Berkeley in 1999, and has been a professor at the University of Georgia since 2011. He has thought about topology and taught mathematics at many levels around the world, with postdoctoral positions at the University of Arizona, the Nankai Institute of Mathematics, and the University of Quebec, a Senior Lectureship at the University of Cape Town, and, most recently, a visiting position as Hirzebruch Research Chair at the Max Planck Institute of Mathematics in Bonn in 2019-20.

\bibliographystyle{plain}
%
%
\bibliography{KirbyCelebratio}

\end{document}